\documentclass[11pt, a4paper]{article}
\usepackage[english]{babel}
\usepackage{cite}
\usepackage{graphics}
\usepackage{graphicx}
\usepackage{epsfig}
\usepackage{amssymb}
\usepackage{hyperref}
\usepackage{authblk}
\usepackage{amsmath}
\usepackage{amsthm}
\usepackage{amsfonts}
\usepackage{float}
\usepackage{booktabs}
\usepackage{enumitem}
\newtheorem{theorem}{Theorem}[section]

 \newtheorem{example}[theorem]{Example}
  \newtheorem{definition}[theorem]{Definition}

\usepackage{siunitx}
\newcommand{\be}{\begin{equation}}
\newcommand{\ee}{\end{equation}}
\newcommand{\ba}{\begin{array}}
\newcommand{\ea}{\end{array}}

\newcommand{\al}{{\alpha}}

\newcommand{\G}{{\Gamma}}

\newcommand{\comment}[1]{}

\newcommand{\La}{\Lambda}

\newcommand{\sg}{\sigma}
\newcommand{\tl}{\tilde}

 \newcommand{\B}[1]{\mbox{\boldmath $#1$}}

\newcommand{\diag}{\mbox{\rm {diag}}}

\newcommand{\vecc}{\mbox{\rm {vec}}}

\begin{document}

\title{Efficient Solution of Parameter Dependent Quasiseparable  Systems
and Computation of Meromorphic Matrix Functions}
\author[a,b]{P.~Boito}

\author[c]{Y.~Eidelman}

\author[d]{L.~Gemignani}

\affil[a]{{\em XLIM--MATHIS, UMR 7252 CNRS Universit\'e de Limoges, 123 
avenue Albert Thomas, 87060 
 Limoges Cedex, France. Email:} paola.boito@unilim.fr}
 
\affil[b]{{\em CNRS, Universit\'e de Lyon, Laboratoire LIP 
(CNRS, ENS Lyon, Inria, UCBL), 46 all\'ee d'Italie, 69364 Lyon Cedex 07, 
France.}}
 
\affil[c]{{\em School of Mathematical Sciences, Raymond and Beverly
Sackler Faculty of Exact Sciences, Tel-Aviv University, Ramat-Aviv,
69978, Israel. Email:} eideyu@post.tau.ac.il}

\affil[d]{{\em Dipartimento di Informatica, Universit\`{a} di Pisa,
Largo Bruno Pontecorvo, 3 - 56127 Pisa, Italy. Email:} l.gemignani@di.unipi.it}

\date{}

\maketitle

\begin{abstract}
In this paper we focus  on the solution of shifted quasiseparable systems  
and of more general parameter dependent matrix equations with quasiseparable 
representations.  We propose an efficient algorithm exploiting the invariance 
of the quasiseparable structure under  diagonal shifting and inversion. This
algorithm is applied to compute various functions of matrices. Numerical 
experiments  show that this approach is fast and numerically robust.  
 
\vspace{5pt}



\end{abstract}

\section{Introduction}
In this paper we propose a novel method for  computing the solution of shifted 
quasiseparable systems and of more general  parameter dependent  linear  matrix
equations with quasiseparable representations. We show that our approach has  
also a noticeable potential for effectively solving some large-scale  algebraic
problems that reduce to evaluating the action of a quasiseparable matrix 
function to a vector. 

Quasiseparable matrices, characterized by the property that off-diagonal blocks have low rank, have found their application in  several  branches of 
applied mathematics and engineering. In the last decade there has been major 
interest in developing fast algorithms for working with such matrices
  \cite{book_egh1,book_egh2,book_vvm1,book_vvm2}. The 
quasiseparable structure arises in the discretization of
 continuous operators, as a consequence of the  local properties of the  
discretization schemes,   and/or because of the decaying
properties of the operator  or its finite approximations. 

It is a celebrated fact that operations with quasiseparable matrices can be 
performed in linear time with respect to their size. In particular,  the QR 
factorization algorithm presented in \cite{eg_qr} computes in linear time a QR
decomposition of a quasiseparable matrix $A\in \mathbb C^{N\times N}$. This decomposition takes the special form $A= V\cdot U\cdot R$, where $R$ is upper triangular,  whereas $U$ and $V$ 
are banded unitary matrices.  It turns out that the matrix $V$ only depends on the generators of the 
strictly lower triangular part of $A$.   This implies that any shifted linear 
system $A+\sigma I_N$, $\sigma\in \mathbb C$, can also be factored as 
$A+\sigma I_N=V\cdot U_\sigma \cdot R_\sigma$  for suitable $U_\sigma$ and 
$R_\sigma$.

Relying upon this fact, in this paper we  design an  efficient algorithm  for  
solving a sequence of shifted quasiseparable linear systems. The invariance of 
the factor $V$  can be exploited to halve the overall computational cost.
Section 2   illustrates the  potential of our approach using motivating 
examples and applications. In Section 3 we describe the novel algorithm by 
proving its correctness. Finally, in Section 4 we present several numerical
experiments where our method is applied to the computation of $f(A)v$ and to the solution 
of linear matrix equations. The results turn out to be as accurate as classical approaches and timings are consistent with the improved complexity estimates.

\section{Motivating Examples}

In this section we describe two motivating examples from applied fields  that 
lead to  the solution of  several shifted or parameter dependent  
quasiseparable linear systems.

\subsection{A Model Problem for Boundary Value ODEs}

Consider the non-local boundary value problem in a linear  finite dimensional 
normed space $X=\mathbb R^N$:
\be\label{eq1}
\frac{d{\B v}}{dt}=A\B v,\quad 0<t< \tau,
\end{equation}
\be\label{eq2}
\frac1{\tau}\int_0^{\tau}\B v(t)\;dt=\B g
\end{equation}
where $A$ is a linear operator in $\mathbb R^N$ and $\B g\in \mathbb R^N$ is a 
given vector.  

The nonlocal problem of the form \eqref{eq1}, \eqref{eq2} has been a subject of  intensive study, see the papers by I.~V.~Tikhonov \cite{T1,T2}, as well as  \cite{ETS} and the literature cited therein.

Under the assumption that all the numbers
$\mu_k=2\pi ik/\tau,\;k=\pm1,\pm2,\pm3,\dots$ are regular points of the 
operator $A$  the problem (\ref{eq1}), (\ref{eq2}) has a unique solution.
Moreover this solution is given by the formula
\begin{equation}\label{eq3}
\B v(t)=q_t(A)\B g, \quad 
q_t(z)=\frac{\tau z e^{zt}}{e^{z\tau}-1}.
\end{equation}
Without loss of generality one can assume that $\tau=2\pi$.

Expanding the function $q_t(z)$ in the Fourier series of $t$ we obtain 
\[
q_t(z)=1+\sum_{k\in{\mathbb Z}/\{0\}}
\frac{ze^{ikt}}{z-ik},\quad 0<t<2\pi.
\]
We consider the equivalent representation with the real series  given by 
\[
q_t(z)=1+2\sum_{k=1}^{\infty}z(z\cos kt-k\sin kt)(z^2+k^2)^{-1},\quad 0<t<2\pi.
\]
Using the formula
\[
-\frac{k}{z^2+k^2}=\frac{z^2}{k(k^2+z^2)}-\frac1{k}
\]
we find that 
\[
q_t(z)=1+2\sum_{k=1}^{\infty}z(z\cos kt+\frac1{k}z^2\sin kt)(z^2+k^2)^{-1}-
2z\sum_{k=1}^{\infty}\frac1{k}\sin kt, \quad 0<t<2\pi
\]
and since
\[
2\sum_{k=1}^{\infty}\frac1{k}\sin kt=\pi-t, \quad 0<t<2\pi
\]
we  arrive at the following formula 
\[
q_t(z)=1+2\sum_{k=1}^{\infty}(z\cos kt+\frac1{k}z^2\sin kt)(z^2+k^2)^{-1}-
z(\pi-t), \quad 0<t<2\pi.
\]

Thus we obtain the  sought expansion  of  the solution $\B v(t)$ of the 
problem \eqref{eq1}, 
\eqref{eq2}:
\be\label{el2}
v(t)=g-(\pi-t)Ag+
2\sum_{k=1}^{\infty}(A\cos kt+\frac1{k}A^2\sin kt)(A^2+  k^2 I_N)^{-1}Ag, 
\quad 0\le t\le2\pi.
\end{equation}
Under the assumptions 
\be\label{ir1}
\|(A-ik I_N)^{-1}\|\le\frac{C}{|k|},\quad k=\pm1,\pm2,\dots
\end{equation} 
and by using the Abel transform 
one can check
that the series in (\ref{el2}) converges uniformly in $t\in[0,2\pi]$. 
Hence,  by  continuity we extend here the formula  \eqref{el2} for the 
solution on all the segment $[0,2\pi]$. The rate of convergence of the 
series in (\ref{el2}) is the same as for the series 
$\sum_{k=1}^{\infty}\frac1{k^2}$. The last may be improved using standard 
techniques for series acceleration  but by including higher degrees of $A$.
Indeed set
$$
C_k(t)=A\cos kt+\frac1{k}A^2\sin kt.
$$
Using the identity
$$
(A^2+k^2 I)^{-1}=\frac1{k^2}I-\frac1{k^2}A^2(A^2+  k^2 I)^{-1}
$$
we get
$$
\sum_{k=1}^{\infty}C_k(t)(A^2+k^2 I_N)^{-1}Ag= 
\sum_{k=1}^{\infty}\frac1{k^2}C_k(t)Ag-\sum_{k=1}^{\infty}\frac1{k^2}C_k(t)
(A^2+k^2 I)^{-1}A^3g.
$$
The first entry here has the form
$$
\sum_{k=1}^{\infty}\frac1{k^2}C_k(t)Ag=
\sum_{k=1}^{\infty}\left(\frac{\cos kt}{k^2}A^2g+
\frac{\sin kt}{k^3}A^3g\right).
$$
Using the formulas
$$
\sum_{k=1}^{\infty}\frac{\cos kt}{k^2}=\frac{\pi^2}{6}-
\frac{\pi}{2}t+\frac{t^2}{4},\quad
\sum_{k=1}^{\infty}\frac{\sin kt}{k^3}=\frac{\pi^2}{6}t-
\frac{\pi}{4}t^2+\frac{t^3}{12},
$$
we get
\be\label{ellu}
v(t)=V_0(t)g+V_1(t)Ag+V_2(t)A^2g+V_3(t)A^3g-
2\sum_{k=1}^{\infty}\frac1{k^2}C_k(t)(A^2+k^2 I)^{-1}A^3g.
\end{equation}
with 
$$
V_0(t)=1,\;V_1(t)=t-\pi,\;
V_2(t)=\frac{\pi^2}{3}-{\pi}t+\frac{t^2}{2},\;
V_3(t)=\frac{\pi^2}{3}t-\frac{\pi}{2}t^2+\frac{t^3}{6}.
$$

Summing up,  our proposal consists in  approximating  the  solution 
$\B v(t)$ of the problem \eqref{eq1}, 
\eqref{eq2}  by the finite sum 
\be
v_\ell(t)=g-(\pi-t)Ag+
2\sum_{k=1}^{\ell}(A\cos kt+\frac1{k}A^2\sin kt)(A^2+k^2 I_N)^{-1}Ag, 
\quad 0\le t\le 2\pi, \label{first_approx}
\end{equation}
or using (\ref{ellu}) by the sum
\begin{gather}
\hat{v}_\ell(t)=\sum_{j=0}^3V_j(t)A^jg \nonumber\\
-2\sum_{k=1}^{\ell}\frac1{k^2}(A\cos kt+\frac1{k}A^2\sin kt)
(A^2+k^2 I_N)^{-1}A^3g,\quad 0\le t\le 2\pi, \label{second_approx}
\end{gather}
where $\ell$ is suitably  chosen by checking the convergence of the expansion.

The computation of
$\B v_\ell(t_i)\simeq \B v(t_i)$, $0\leq i\leq M+1$,  requires the solution  
 of a  possibly large set of the shifted systems
of the form
\be\label{len}
(A+\sg_iI_N)\B x_i=\B y,\quad i=1,\dots,\ell.
\end{equation}

The same conclusion applies  to the problem of computing the function of  a 
quasiseparable matrix whenever  the function can be represented as a series of 
partial fractions. The classes of meromorphic functions admitting such a 
representation were investigated for instance in \cite{SH}. Other 
partial fraction approximations of
certain analytic functions can be found in \cite{HHT}.
In the next section we describe an effective algorithm for this task. 

As additional context for this model problem, note that a numerical approximation of the solution can be obtained by using the 
discretization of the boundary value problem on a grid and the subsequent 
application of the cyclic reduction approach discussed in \cite{AP}. 
This approach can be combined with techniques for preserving  an approximate  quasiseparable structure in 
recursive LU-based solvers: see the recent papers
\cite{BK,BMR,GZ}. The approximate 
quasiseparable structure of functions of quasiseparable matrices has also been 
investigated in \cite{massei2016decay}.

\subsection{Sylvester-type  Matrix Equations}\label{subsec_matrix_equations}
As a natural extension of the problem \eqref{len}, the right-hand side $\B y$ 
could also depend on the parameter $\sigma$,
that is, $\B y=\B y(\sigma)$ and $\B y_i=\B y(\sigma_i)$, $i=1,\dots,\ell$.  
This situation is
common in many applications  such as control theory, structural dynamics and time-dependent PDEs \cite{GS}.
In this case, the systems to be solved take the form
\[
A X +X D =Y, \quad A\in \mathbb R^{N\times N}, \ D=\diag\left[\sigma_1, \ldots, \sigma_\ell\right],
\ Y=\left[\B y_1, \ldots, \B y_\ell \right].
\]
Using the Kronecker product  this matrix equation can be rewritten as a bigger linear system 
$\mathcal A\, \vecc(X)=\vecc(Y)$, where  $\mathcal A=I_\ell \otimes A + D^T \otimes I_N \in \mathbb R^{N \ell \times N \ell}$,
$\vecc(X)=\left[\B x_1^T, \ldots, \B x_\ell^T\right]^T$, $\vecc(Y)=\left[\B y_1^T, \ldots, \B y_\ell^T\right]^T$.

The extension to the case where $D$ is replaced by a lower triangular matrix 
$L=(l_{i,j})\in\mathbb R^{\ell \times \ell}$ is 
based on the backward substitution technique  which amounts to solve
\be\label{len1}
(A+l_{i,i} I_N)\B x_i =\B y_i -\sum_{j=i+1}^\ell l_{i,j} \B x_j, \quad i=\ell \colon -1 \colon 1.
\end{equation} 
Such approach  has been used in different related contexts  where the considered Sylvester equation is
occasionally called a {\em sparse-dense} equation \cite{S}.  The classical reduction proposed by
Bartels and Stewart \cite{BS} makes it possible to  deal with a general matrix $F$ by  first computing
 its Schur decomposition $F=U L U^H$ and then solving  $A (X U) + (XU) L = Y U$.  The resulting approach is well
suited especially  when  the size of $A$ is large w.r.t. the number of shifts.
 If $A$ is  quasiseparable  then \eqref{len1}  again reduces to computing a 
sequence of shifted  systems having the same
structure in the lower triangular part and the method presented in the next 
section can be used.

\section{The basic algorithm}\label{sec_algo}

Let us first recall the definition of quasiseparable matrix structure and quasiseparable generators. See \cite{book_egh1} for more details.

\begin{definition}
A block matrix  $A=(A_{i,j})_{i,j=1}^N$, with block entries
$A_{i,j}\in \mathbb R^{m_i\times m_j}$, is said to have {\em lower quasiseparable generators}
$P(i) \in \mathbb R^{m_i\times r^L_{i-1}} \ (i=2,\dots,N)$, $Q(j)\in \mathbb R^{r^L_j \times m_j} \ (j=1,\dots,N-1)$,
$\Xi(k)\in \mathbb R^{r^L_{k}\times r^L_{k-1}} \ (k=2,\dots,N-1)$
of orders $r^L_k\;(k=1,\dots,N-1)$ and  {\em upper quasiseparable generators}
$G(i)\in \mathbb R^{m_i\times r^U_i} (i=1,\dots,N-1)$, 
$H(j)\in \mathbb R^{r^U_{j-1} \times m_j}\ (j=2,\dots,N)$, $\Theta(k)\in \mathbb R^{r^U_{k-1}\times r^U_k} \,(k=2,\dots,N-1)$ of orders
$r^U_k\;(k=1,\dots,N-1)$  if
\[
A_{i,j}=\left \{ \begin{array}{ll}
  P(i)\Xi_{i,j}^>Q(j) \ {\rm if}\ 1\leq j<i\leq N,\\
  G(i)\Theta_{i,j}^<H(j) \ {\rm if}\ 1\leq i<j\leq N,
\end{array}\right.
\]
where $\Xi_{i,j}^>=\Xi(i-1)\cdots \Xi(j+1)$ for $i>j+1$ and $\Xi_{j+1,j}=I_{r^L_j}$, and,
similarly, $\Theta_{i,j}^<=\Theta(i+1)\cdots \Theta(j-1)$  for $j>i+1$ and  $\Theta_{i,i+1}=I_{r^U_i}$.

If $A$ admits such a representation, is is said to be {\em $(r_L,r_U)$-quasiseparable}.
Diagonal entries are stored separately, that is, we set $\Lambda(i)=A_{i,i}\in\mathbb{R}^{m_i\times m_i}$. 

The same representation can be applied to complex matrices.

\end{definition}

We have denoted the quasiseparable generators by capital letters, as it is often done for matrices. Note however that the generators can be numbers, vectors or matrices, depending on the quasiseparable orders $r_{L_i}, r_{U_j}$  and on the block sizes $m_i, m_j$.

To solve the systems \eqref{len},  \eqref{len1} we rely upon the
QR factorization algorithm described in \cite{eg_qr} (see also Chapter 20
of \cite{book_egh1}). At first we compute  the factorization 
\be\label{ts}
A+\sg I=V\cdot T_{\sg}
\end{equation}
with a unitary matrix $V$ and a lower banded, or a block upper triangular, 
matrix $T_{\sg}$. It turns out that the matrix $V$ does not depend on 
$\sg$ at all
and moreover an essential part of the quasiseparable generators of the matrix $T_\sg$ do
not depend on $\sg$ either. So a relevant  part of the computations for all the 
values of $\sg$ only needs to be performed once. Thus the problem is reduced to 
the solution of the set of the systems
\be\label{len11}
T_{\sg}\B x_{\sg}=V^H\B y_i,\quad \sg=\sg_i,\;i=1,\dots,\ell.
\end{equation}
The inversion of every matrix $T_{\sg}$ as well as the solution of the 
corresponding linear system is significantly simpler than for the original 
matrix. We compute the factorization
\be\label{eqTUR}
T_{\sg}=U_{\sg}R_{\sg}
\end{equation}
with a block upper triangular unitary matrix $U_{\sg}$ and upper triangular 
$R_{\sg}$ and solve the systems
\be\label{eqtriang}
R_{\sg}\B x_{\sg}=U_{\sg}^HV^H\B y_i.
\end{equation}
Thus we obtain our main algorithm, which takes as input a set of quasiseparable generators for $A$, shifts $\sigma_i$ and right-hand block vectors $\B y_i$, and outputs the solutions $\B x_i$ of the linear systems $(A+\sg_i I)\B x_i=\B y_i$. The algorithm is comprised of two parts.
\begin{itemize}
\item Part 1 computes useful quantities that are common to all the linear systems, namely, quasiseparable generators for $V$ in the factorization \eqref{ts}, and some quasiseparable generators for each $T_{\sg_i}$ that do not actually depend on $\sigma_i$. 
\item Part 2 uses the results of Part 1, along with input data, to solve efficiently each linear system (which is why it begins with a loop over all the systems). For each $i=1,\ldots,\ell$ it computes the factorization \eqref{eqTUR}, and then solves the triangular system \eqref{eqtriang}.
\end{itemize}
Note that each of the matrices $V$ and $U_{\sg_i}$ can also be factored as the product of $N$ ``small'' unitary matrices, which are denoted as $V_k$ and $U_k^{(i)}$, respectively, in the algorithm that follows. Throughout the algorithm, these factored representations are computed via successive QR factorizations of suitable matrices and then used to compute products by $V^H$ or $U_{\sg_i}^H$ in $O(N)$ time: see the proof of the algorithm for more details.

The sentences in italics explain the purpose of each block of instructions.


\medskip

\begin{center}
{\bf Main algorithm}
\end{center}

\medskip

{\bf Input:} 
\begin{itemize}
\item quasiseparable generators for the block matrix $A$, with entries $(A_{i,j})_{i,j=1}^N$ of sizes $m_i\times m_j$, namely:
\begin{itemize}
\item lower quasiseparable generators
$P(i)\;(i=2,\dots,N),\;Q(j)\;(j=1,\dots,N-1),\;\Xi(k)\;(k=2,\dots,N-1)$
of orders $r^L_k\;(k=1,\dots,N-1)$,
\item upper quasiseparable generators
$G(i)\; (i=1,\dots,N-1),\;
H(j)\;(j=2,\dots,N),\;\Theta(k)\;(k=2,\dots,N-1)$ of orders
$r^U_k\;(k=1,\dots,N-1)$,
\item diagonal entries
$\Lambda(k)\;(k=1,\dots,N)$;
\end{itemize}
\item complex numbers $\sg_i,\;i=1,\dots,\ell$;
\item block vectors $\B y_i=(y_i(k))_{k=1}^N$ with $m_k$-
dimensional coordinates $y(k)$.
\end{itemize}

{\bf Output:} solutions $\B x^{(i)}=\B x_{\sigma_i},\; i=1,\dots,\ell$ of the systems \eqref{len11}.

\medskip

{\bf Part 1.}

\begin{itemize}

\item[1.] {\em Initialize auxiliary quantities:}
\[
\rho_N=0,\quad\rho_{k-1}=\min\{m_k+\rho_k,\;r^L_{k-1}\},\quad k=N,\dots,2,\;
\rho_0=0,
\]
\[
\rho'_k=\rho_k+r^U_k,\;k=1,\dots,N-1,\quad    
\nu_k=m_k+\rho_k-\rho_{k-1},\;k=1,\dots,N.
\]

\item[2.] {\em Compute the quasiseparable representation of the matrix $V$ (that is, generators with subscript $_V$) and the basic 
elements of the representation of the matrix $T_{\sg}$ (that is, generators with subscript $_T$), as well as the vector
$\B w_i=V^H\B y_i$. Note that this is done through the computation of the ``small'' unitary factors $V_k$ of $V$.}

\begin{itemize}
\item Using the QR factorization or another algorithm compute the factorization
\be\label{aln}
P(N)=V_N\left(\ba{c}X_N\\0_{\nu_N\times r^L_{N-1}}\end{array}\right),
\end{equation}
where $V_N$ is a unitary matrix of sizes $m_N\times m_N$, $X_N$ is a
matrix of sizes $\rho_{N-1}\times r^L_{N-1}$. 

Determine the matrices $P_V(N),
\Lambda_V(N)$ of sizes $m_N\times\rho_{N-1},m_N\times\nu_N$  from the partition
\be\label{eln}
V_N=\left(\ba{cc}P_V(N)&\La_V(N)\end{array}\right).
\end{equation}

Compute 
\be\label{hgn}
\left(\ba{c}H_T(N)\\ \Lambda_T(N)\ea\right)=\left(\ba{c}H(N)\\V_N^H \La(N)\ea\right),
\end{equation}
\be\label{cn}
\left(\ba{c}c_{N,i}\\w_i(N)\ea\right)=V^H_Ny_i(N), \quad 1\leq i\leq \ell
\end{equation}
with the matrices $H_T(N),\La_T(N),c_{N,i},w_i(N)$ of sizes $\rho'_{N-1}\times m_N,
\nu_N\times m_N,\rho_{N-1}\times1,\nu_N\times1$ respectively.

Set
\be\label{gn}
\Gamma_N=\left(\ba{c}0_{r^U_{N-1}\times m_N}\\P_V^H(N)\ea\right).
\end{equation}

\item For $k=N-1,\dots,2$ perform the following.

Using the QR factorization or another algorithm
compute the factorization
\be\label{al}
\left(\ba{c}P(k)\\X_{k+1}\Xi(k)\end{array}\right)=
V_k\left(\ba{c}X_k\\0_{\nu_k\times r^L_{k-1}}\end{array}\right),
\end{equation}
where $V_k$ is a unitary matrix of sizes $(m_k+\rho_k)\times(m_k+\rho_k)$,
$X_k$ is a matrix of sizes $\rho_{k-1}\times r^L_{k-1}$.

Determine the matrices $P_V(k),Q_V(k),\Xi_V(k),\La_V(k)$ of sizes
$m_k\times\rho_{k-1},\rho_k\times\nu_k,\rho_k\times\rho_{k-1}$,
$m_k\times\nu_k$ from the partition
\be\label{el}
V_k=\left(\ba{cc}P_V(k)&\La_V(k)\\ \Xi_V(k)&Q_V(k)\end{array}\right).
\end{equation}

Compute
\be\label{hg}
\left(\ba{cc}H_T(k)&\Theta_T(k)\\ \La_T(k)&G_T(k)\ea\right)=
\left(\ba{cc}I_{r^U_{k-1}}&0\\0&V_k^H\ea\right)
\left(\ba{ccc}H(k)&\Theta(k)&0\\ \La(k)&G(k)&0\\X_{k+1}Q(k)&0&I_{\rho_k}\ea\right).
\end{equation}
with the matrices $H_T(k),\Theta_T(k),\La_T(k),G_T(k)$ of sizes 
$\rho'_{k-1}\times m_k,\rho'_{k-1}\times\rho'_k,\nu_k\times m_k,
\nu_k\times\rho'_k$ respectively.

Set
\be\label{g}
\Gamma_k=\left(\ba{c}0_{r^U_{k-1}\times m_k}\\P_V^H(k)\ea\right)
\end{equation}
and compute
\be\label{c}
\left(\ba{c}c_{k,i}\\w_i(k)\ea\right)=V_k^H\left(\ba{c}y_i(k)\\c_{k+1,i}\ea\right), \quad 1\leq i\leq \ell
\end{equation}
with the vector columns $c_{k,i},w_i(k)$ of sizes $\rho_{k-1},\nu_k$ respectively.

\item Set $V_1=I_{\nu_1}$ and 
\be\label{hg1}
\La_T(1)=\left(\ba{c}\La(1)\\X_2 Q(1)\ea\right),\;
G_T(1)=\left(\ba{cc}G(1)&0\\0&I_{\rho_1}\ea\right),\quad
\Gamma_1=\left(\ba{c}I_{m_1}\\0_{\rho_1\times m_1}\ea\right),
\end{equation}
\be\label{c1}
w_i(1)=V_1^H\left(\ba{c}y_i(1)\\c_{2,i}\ea\right), \quad 1\leq i\leq \ell.
\end{equation}
\end{itemize}
\end{itemize}

\noindent

{\bf Part 2.}

For $i=1,\dots,\ell$ {\em (that is, for each shifted linear system)} perform the following:

\begin{itemize}
\item[3.] {\em Compute the factorization $T_{\sg_i}=U_{\sg_i}R_{\sg_i}$ and the vector
$\B v^{(i)}=\B v_{\sg_i}=U_{\sg_i}^H\B w_i$, $\B w_i=V^H \B y_i$. In particular, compute the ``small'' unitary factors $U_k^{(i)}$ of $U_{\sg_i}$ and use this factorization to determine the quasiseparable generators of $T_{\sg_i}$, denoted by the subscript $_T$, and the vector $\B v^{(i)}$.}
\begin{itemize}
\item 
Compute the QR factorization
\be\label{d1}
\La_T(1)+\sg_i\Gamma_1=U_1^{(i)}
\left(\ba{c}\La_R^{(i)}(1)\\0_{\rho_1\times m_1}\ea\right),
\end{equation}
where $U_1^{(i)}$ is a $\nu_1\times\nu_1$ unitary matrix and $\La_R^{(i)}(1)$ is 
an upper triangular $m_1\times m_1$ matrix.

Compute
\be\label{g1}
\left(\ba{c}G_R^{(i)}(1)\\Y_1^{(i)}\ea\right)=(U_1^{(i)})^H G_T(1),
\end{equation}
\be\label{vv1}
\left(\ba{c}v^{(i)}(1)\\\al^{(i)}_1\ea\right)=(U_1^{(i)})^Hw_i(1)
\end{equation}
with the matrices $G_R^{(i)}(1),v^{(i)}(1),Y_1^{(i)},\al^{(i)}_1$ of sizes
$m_1\times\rho'_1,m_1\times1,\rho_1\times\rho'_1,\rho_1\times1$.

\item
For $k=2,\dots,N-1$ perform the following.

Compute the QR factorization
\be\label{d}
\left(\ba{c}Y^{(i)}_{k-1}(H_T(k)+\sg_i\Gamma_k)\\ \La_T(k)+\sg_i \La_V^H(k)\ea\right)=
U_k^{(i)}\left(\ba{c}\La_R^{(i)}(k)\\0_{\rho_k\times m_k}\ea\right),
\end{equation}
where $U_k^{(i)}$ is an $(m_k+\rho_k)\times(m_k+\rho_k)$ unitary matrix and
$\La_R^{(i)}(k)$ is an $m_k\times m_k$ upper triangular matrix.

Compute
\be\label{gg}
\left(\ba{c}G_R^{(i)}(k)\\Y_k^{(i)}\ea\right)=
(U_k^{(i)})^H\left(\ba{c}Y^{(i)}_{k-1}\Theta_T(k)\\G_T(k)\ea\right),
\end{equation}
\be\label{vv}
\left(\ba{c}v^{(i)}(k)\\\al^{(i)}_k\ea\right)=
(U_k^{(i)})^H\left(\ba{c}\al^{(i)}_{k-1}\\w_i(k)\ea\right)
\end{equation}
with the matrices $G_R^{(i)}(k),v^{(i)}(k),Y_k^{(i)},\al^{(i)}_k$ of sizes
$m_k\times\rho'_k,m_k\times1,\rho_k\times\rho'_k,\rho_k\times1$.

\item
Compute the $QR$ factorization
\be\label{dn}
\left(\ba{c}Y^{(i)}_N H_T(N)+\sg_i\Gamma_N\\ \La_T(N)+\sg_i \La_V^H(N)\ea\right)=
U^{(i)}_N \La_R^{(i)}(N),
\end{equation}
where $U^{(i)}_N$ is a unitary matrix of sizes
$(\nu_N+\rho_{N-1})\times(\nu_N+\rho_{N-1})$ and $\La_R^{(i)}(N)$ is an upper
triangular matrix of sizes $m_N\times m_N$.

Compute
\be\label{vvn}
v^{(i)}(N)=(U^{(i)}_N)^H\left(\ba{c}\al^{(i)}_{N-1}\\w_i(N)\ea\right).
\end{equation}
\end{itemize}

\item[4.]
 {\em Solve the system $R_{\sg_i}\B x^{(i)}=\B v^{(i)}$, using the previously computed quasiseparable generators.}
\begin{itemize}
\item 
Compute
\[
x^{(i)}(N)=(\La_R^{(i)}(N))^{-1}v^{(i)}(N),\]
\[\eta^{(i)}_{N-1}=(H_T(N)+\sg_i\Gamma_N)v^{(i)}(N)
\]

\item
For $k=N-1,\dots,2$ compute
$$
x^{(i)}(k)=(\La_R^{(i)}(k))^{-1}(v^{(i)}(k)-G_R^{(i)}(k)\eta^{(i)}_k),$$
$$\eta^{(i)}_{k-1}=\Theta_T(k)\eta^{(i)}_k+(H_T(k)+\sg_i\Gamma_k)x^{(i)}(k).
$$

\item
Compute
$$
x^{(i)}(1)=(\La_R^{(i)}(1))^{-1}(v^{(i)}(1)-G_R^{(i)}(1)\eta^{(i)}_1)
$$

\end{itemize}
\end{itemize}

{\em Proof of correctness.} 
The shifted matrix $A+\sg I_N$ has the same lower and upper quasiseparable 
generators as the matrix $A$ and diagonal entries $\La(k)+\sg I_{m_k}\;(k=1,\dots,N)$.
To compute the factorization (\ref{ts}) we apply Theorem 20.5 from \cite{book_egh1}, obtaining the representation of the unitary matrix
$V$ in the form
\be\label{vn}
V=\tl V_N\tl V_{N-1}\cdots\tl V_2\tl V_1,
\end{equation}
where
$$
\tl V_1=V_1\oplus I_{\phi_1},\quad
\tl V_k=I_{\eta_k}\oplus V_k\oplus I_{\phi_k},\;k=2,\dots,N-1,\quad
\tl V_N=I_{\eta_N}\oplus V_N 
$$
with $\eta_k=\sum_{j=1}^{k-1}m_j,\;\phi_k=\sum_{j=k+1}^N m_j$ and 
$(m_k+\rho_k)\times(m_k+\rho_k)$ unitary matrices $V_k$, as well as the 
formulas (\ref{aln}), (\ref{eln}), (\ref{al}), (\ref{el}), $V_1=I_{\nu_1}$.
Here we see that the matrix $V$ does not depend on $\sg$. Moreover the
representation (\ref{vn}) yields the formulas (\ref{cn}), (\ref{c}), (\ref{c1})
to compute the vectors  $\B w_i=V^H\B y_i$, $1\leq i\leq \ell$.

Next we apply the corresponding formulas from the same theorem
to compute diagonal entries $\La_T^{\sg}(k)\;(k=1,\dots,N)$ and upper 
quasiseparable generators $G_T^{\sg}(i)\;(i=1,\dots,N-1),\;
H_T^{\sg}(j)\;(j=2,\dots,N),\;\Theta_T^{\sg}(k)\;(k=2,\dots,N-1)$ of the matrix $T_\sigma$.   Hence, we obtain that 
\begin{gather*}
\left(\ba{c}H_T^{\sg}(N)\\ \La_T^{\sg}(N)\ea\right)=
\left(\ba{cc}I_{r^U_{N-1}}&0\\0&V_N^H\ea\right)
\left(\ba{c}H(N)\\ \La(N)+\sg I_{m_N}\ea\right),\\
\left(\ba{cc}\La_T^{\sg}(1)&G_T^{\sg}(1)\ea\right)=
\left(\ba{ccc}\La(1)+\sg I_{m_1}&G(1)&0\\
X_2 Q(1)&0&I_{\rho_1}\ea\right), 
\end{gather*}
and for $k=N-1, \ldots, 2$, 
\[\left(\ba{cc}H_T^{\sg}(k)& \Theta_T^{\sg}(k)\\ \La_T^{\sg}(k)&G_T^{\sg}(k)\ea\right)=
\left(\ba{cc}I_{r^U_{k-1}}&0\\0&V_k^H\ea\right)
\left(\ba{ccc}H(k)&\Theta(k)&0\\ \La(k)+\sg I_{m_k}&G(k)&0\\
X_{k+1}Q(k)&0&I_{\rho_k}\ea\right).
\]
From here we obtain the formulas
\be\label{aug}
\begin{gathered}
H_T^{\sg}(k)=H_T(k)+\sg\Gamma_k,\; k=N,\dots,2,\\
\La_T^{\sg}(k)=\La_T(k)+\sg \La_V^*(k),\; k=N,\dots,2,\quad 
\La_T^{\sg}(1)=\La_T(1)+\sg\Gamma_1,\\
G_T^{\sg}(k)=G_T(k),\;k=1,\dots,N-1,\quad
\Theta_T^{\sg}(k)=\Theta_T(k),\;k=2,\dots,N-1
\end{gathered}
\end{equation}
with $H_T(k),\La_T(k),G_T(k),\Theta_T(k)$ and $\Gamma_k$ as in (\ref{hgn}), (\ref{gn}),
(\ref{hg}), (\ref{g}) and (\ref{hg1}).

Now  by applying Theorem 20.7 from \cite{book_egh1} to the matrices $T_{\sg_i},\;i=1,2,
\dots,M$ with the generators determined in (\ref{aug}) we obtain the formulas 
(\ref{d1}), (\ref{g1}), (\ref{d}), (\ref{gg}), (\ref{dn}) to compute unitary 
matrices $U_k^{(i)}$ and quasiseparable generators of the upper triangular 
$R_{\sg_i}$ such that $T_{\sg_i}=U_{\sg_i}R_{\sg_i}$, where
\be\label{un}
U_{\sg_i}=\tl U_1^{(i)}\tl U_2^{(i)}\cdots\tl U_{N-1}^{(i)}\tl U_N^{(i)}
\end{equation}
with
$$
\tl U^{(i)}_1=U^{(i)}_1\oplus I_{\phi_1},\quad
\tl U^{(i)}_k=I_{\eta_k}\oplus U^{(i)}_k\oplus I_{\phi_k},\;k=2,\dots,N-1,\quad
\tl U^{(i)}_N=I_{\eta_N}\oplus U^{(i)}_N. 
$$
Moreover the representation (\ref{un}) yields the formulas (\ref{vv1}), 
(\ref{vv}), (\ref{vvn}) to compute the vector $\B v^{(i)}=U^H_{\sg_i}\B w_i$, $1\leq i\leq l$.

Finally applying Theorem 13.13 from \cite{book_egh1} we obtain Step 2.2 to compute 
the solutions of the systems $R_{\sg_i}\B x^{(i)}=\B v_{(i)}$.
$\hfill\Box$

Concerning the complexity of the previous algorithm  we observe that,   under the
 simplified assumptions of  $r^L_k=r^U_k=r$,  $m_i=m_j=m$ and $r\ll m$, the   cost of step 1  is of the order
$(6r^2m + 2m^2)(Nm)$,  whereas the cost of step 2 can be estimated as  $(2m^2\ell)(Nm)$ arithmetic operations.
Therefore, the proposed algorithm saves at least half of the overall cost of solving $\ell$ shifted
quasiseparable linear systems.

\section{Numerical Experiments}

The proposed fast algorithm has been implemented in MATLAB.\footnote{The code is available at \href{http://www.unilim.fr/pages_perso/paola.boito/software.html}{\tt http://www.unilim.fr/pages\_perso/paola.boito/software.html}.} All the experiments were performed on a MacBookPro equipped with MATLAB R2016b.

\begin{figure}
\centering
\includegraphics[width=0.9\textwidth]{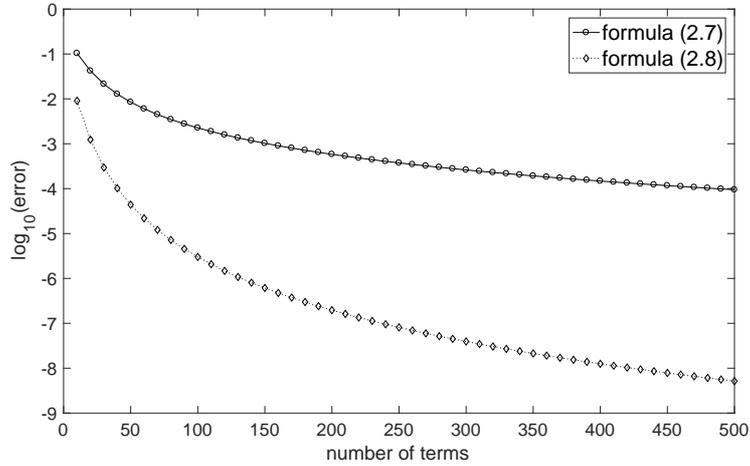}
\caption{Results for Example \ref{ex_series} for $t=\pi/2$. This is a logarithmic plot of the errors given by the application of formulas \eqref{first_approx} (circles with solid line) and \eqref{second_approx} (diamonds with dotted line) for the computation of a matrix function times a vector.}\label{fig_ex1_errors_pi2} 
\end{figure}

\begin{figure}
\centering
\includegraphics[width=0.9\textwidth]{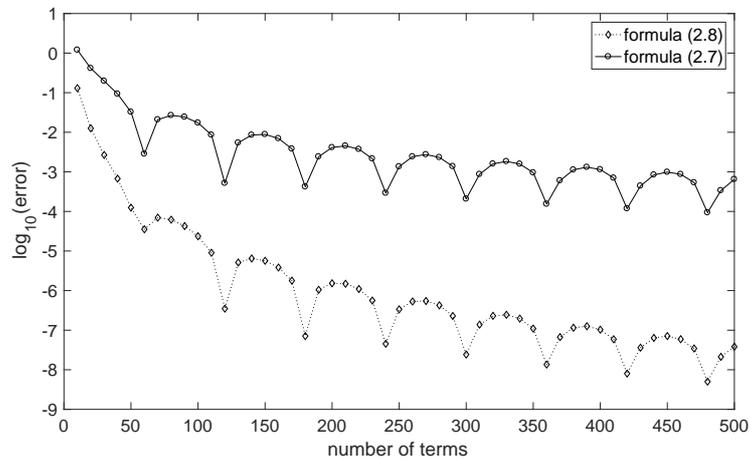}
\caption{Results for Example \ref{ex_series} for $t=\pi/12$. This is a logarithmic plot of the errors given by the application of formulas \eqref{first_approx} (circles with solid line) and \eqref{second_approx} (diamonds with dotted line) for the computation of a matrix function times a vector. In this example the errors do not decrease monotonically.}\label{fig_ex1_errors_pi12} 
\end{figure}

\begin{example}\label{ex_series}
Let us test the computation of functions of quasiseparable matrices via series expansion, as outlined in Section 2. We choose $A$ as the $100\times 100$ one-dimensional discretized Laplacian with zero boundary conditions, which is $(1,1)$-quasiseparable, and $g$ as a random vector with entries taken from a uniform distribution over $[0,1]$. Define 
$$
v_{ex}=2\pi A e^{At}(e^{2\pi A}-I)^{-1}g,
$$ 
as exact solution (computed in multiprecision) to the problem \eqref{eq1}, \eqref{eq2}. Let $v_{\ell}$ and $\hat{v}_{\ell}$ be the approximate solutions obtained from \eqref{first_approx} and \eqref{second_approx}, respectively, with $\ell$ expansion terms. Figures \ref{fig_ex1_errors_pi2} and \ref{fig_ex1_errors_pi12} show logarithmic plots of the absolute normwise errors $\|v_{ex}-v_{\ell}\|_2$ and $\|v_{ex}-\hat{v}_{\ell}\|_2$ for $t=\pi/2$ and $t=\pi/12$, and values of $\ell$ ranging from $10$ to $500$. The results clearly confirm that the formulation \eqref{second_approx} has  improved convergence properties with respect to \eqref{first_approx}. Note that the decreasing behavior of the errors is not always monotone. 

It should be pointed out that, in this approach, a faster convergence of the series expansion gives a faster method for approximating the solution vector with a given accuracy. Indeed, the main computational effort comes from the solution of the shifted linear systems $(A^2+k^2 I_N)\B x_k =Ag$ or $(A^2+k^2 I_N)\B x_k =A^3g$, and it is therefore proportional to the number of terms in the truncated expansion.
\end{example}

\begin{example}\label{ex_simoncini}
This example is taken from \cite{simoncini2010extended}, Example 3.3. We consider here the matrix $A\in\mathbb{R}^{2500\times 2500}$ stemming from the centered finite difference discretization of the differential operator $-\Delta u+10u_x$ on the unit square with homogeneous Dirichlet boundary condition. Note that $A$ has (scalar) quasiseparability order $50$, but it can also be seen as block $1$-quasiseparable with block size $50$. 

We want to compute the matrix function
$$ A^{\frac12}b\approx \sum_{k=1}^{\ell} \kappa_k(\omega_k^2I-A)^{-1}Ab,$$ 
where $b$ is the vector of all ones and the choice of the coefficients $\kappa_k$ and $\omega_k$ corresponds to the choice of a particular rational approximation for the square root function. 

We apply the rational approximations presented in \cite{HHT} as {\tt method2} and {\tt method3}; the latter is designed specifically for the square root function and it is known to have better convergence properties. Just as for Example \ref{ex_series}, the main computational burden here consists in computing the terms $(\omega_k^2I-A)^{-1}Ab$ for $k=1,\ldots,\ell$, that is, solving the $\ell$ shifted quasiseparable linear systems $(\omega_k^2I-A)x_k=Ab$. 

In this example we are especially interested in testing the numerical robustness of our fast solver when compared to classical solvers available in MATLAB, in the context of computing matrix functions.  

Table \ref{table_simoncini} shows $2$-norm relative errors with respect to the result computed by the MATLAB command {\tt sqrtm(A)*b}, for several values of $\ell$ (number of terms in the expansion or number of integration nodes). We have tested three approaches to solve the shifted linear systems involved in the expansion: classical backslash solver, fast structured solver with blocks of size $1$ and quasiseparability order $50$, and fast structured solver with blocks of size $50$ and quasiseparability order $1$. For each value of $\ell$, the errors are roughly the same for all three approaches (so a single error is reported in the table). In particular, the fast algorithms appear to be as accurate as standard solvers. 

\end{example}

\begin{table}
\caption{Relative errors for Example \ref{ex_simoncini}.\label{table_simoncini}}
\centering
\begin{tabular}{|c|c|c|c|}
\toprule
 $\ell$ & rel. error ({\tt method 2}) & rel. error ({\tt method 3})\\
\midrule
$6$ & $4.58$e$-4$ & $3.25$e$-5$ \\
$8$ & $3.03$e$-5$ & $6.77$e$-7$\\
$10$ & $9.67$e$-7$ & $9.90$e$-9$\\
$12$ & $1.55$e$-8$ & $1.16$e$-10$\\
$14$ & $2.06$e$-9$ & $1.56$e$-12$\\
$16$ & $2.25$e$-11$ & $2.16$e$-13$\\
$18$ & $2.29$e$-12$ & $2.09$e$-13$\\
$20$ & $2.44$e$-13$ & $2.14$e$-13$\\
\bottomrule
\end{tabular}
\end{table}

\begin{example}\label{ex_poisson}
The motivation for this example comes from the classical problem of solving the Poisson equation on a rectangular domain with uniform zero boundary conditions. The equation takes the form
$$\Delta u(x,y)=f(x,y),\quad \textrm{with}\quad 0<x<a,\, 0<y<b,$$
and its finite difference discretization yields a matrix equation
\begin{equation}\label{rect_poisson}
AX+XB=F\quad \textrm{with}\quad X,F\in\mathbb{R}^{N_b\times N_a}, 
\end{equation}
where $N_a$ is the number of grid points taken along the $x$ direction and  $N_b$ is the number of grid points along the $y$ direction. Here $A$ and $B$ are matrices of sizes $N_b\times N_b$ and $N_a\times N_a$, respectively, and both are Toeplitz symmetric tridiagonal with nonzero entries $\{-1,2,1\}$. See e.g., \cite{wan1973core} for a discussion of this problem.

A widespread approach consists in reformulating the matrix equation \eqref{rect_poisson} as a larger linear system of size $N_aN_b\times N_aN_b$ via Kronecker products. Here instead we apply the idea outlined in Section \ref{subsec_matrix_equations}: compute the (well-known) Schur decomposition of $B$, that is, $B=UDU^H$, and solve the equation $A(XU)+(XU)D=FU$. Note that, since $D$ is diagonal, this equation can be rewritten as a collection of shifted linear systems, where the right-hand side vector may depend on the shift. This approach is especially interesting when $N_b$ is significantly larger than $N_a$.

Table \ref{table_poisson} shows relative errors on the solution matrix $X$, computed w.r.t. the solution given by a standard solver applied to the Kronecker linear system. Here we take $F$ as the matrix of all ones. The results show that our fast method computes the solution with good accuracy.
\end{example}

\begin{table}
\caption{Relative errors for Example \ref{ex_poisson}.\label{table_poisson}}
\centering
\begin{tabular}{|c c| c c c c c|}
\toprule
 &$N_a$ & $10$ & $25$ & $50$ & $75$ & $100$\\ 
$N_b$&&&&&&\\
 \midrule
 $50$ && $9.58$e$-16$ & $1.98$e$-14$ & $2.05$e$-14$ & $9.34$e$-14$& $2.14$e$-13$\\
$100$ && $5.55$e$-15$ & $2.30$e$-14$& $4.58$e$-14$ & $2.11$e$-13$& $5.61$e$-13$\\
$150$ && $4.93$e$-15$ & $3.20$e$-14$& $1.22$e$-13$ & $1.75$e$-13$& $2.36$e$-13$\\
$200$ && $1.25$e$-14$ & $6.48$e$-14$& $2.33$e$-13$ & $3.97$e$-13$& $6.23$e$-13$\\
$250$ && $3.20$e$-15$ & $1.23$e$-14$& $6.80$e$-14$ & $8.98$e$-14$& $1.54$e$-13$\\
$500$ && $3.77$e$-15$ & $1.82$e$-14$& $4.85$e$-14$ & & \\
$1000$ && $6.08$e$-15$ & $3.02$e$-14$&  &&\\
\bottomrule
\end{tabular}
\end{table}

\begin{figure}
\centering
\includegraphics[width=\textwidth]{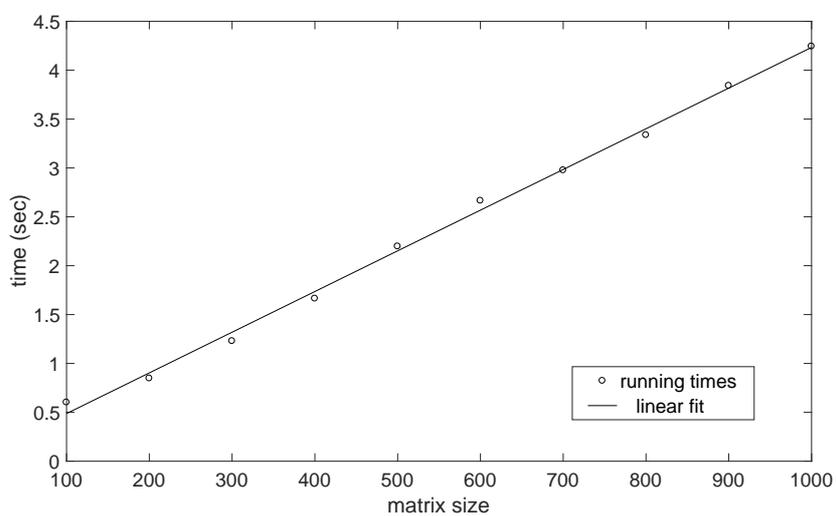}
\caption{Running times for Example \ref{ex_time_size}. The linear fit appears to be a good approximation of the actual data. Its  equation is $y=0.0042x+0.071$.}\label{fig_time_size}
\end{figure}

In the next examples we test experimentally the complexity of our algorithm.
\begin{example}\label{ex_time_size}
We consider matrices $A_n\in\mathbb{R}^{n\times n}$ defined by random quasiseparable generators of order $3$. The second column of Table \ref{table_time_size} shows the running times of our structured algorithm applied to randomly shifted systems $(A_n+\sigma_i I_n)\bf{x}_i=\bf{y}$, for $i=1,\ldots,50$ and growing values of $n$. The same data are plotted in Figure \ref{fig_time_size}: the growth of the running times looks linear with $n$, as predicted by theoretical complexity estimates.

The third column of Table \ref{table_time_size} shows timings for the structured algorithm applied sequentially (i.e., without re-using the factorization \eqref{ts}) to the same set of shifted systems. The gain obtained by the fast algorithm of Section \ref{sec_algo} w.r.t. a sequential structured approach amounts to a factor of about $2$, which is consistent with the discussion at the end of Section \ref{sec_algo}. Experiments with a different number of shifts yield similar results.
\end{example}

\begin{table}
\caption{Running times in seconds for Example \ref{ex_time_size}.\label{table_time_size}}
\centering
\begin{tabular}{|c|c|c|c|}
\toprule
 system size $n$&fast algorithm & sequential algorithm&ratio\\ 
 \midrule
$100$ & $0.6016$ & $1.2049$&$2.0029$\\
$200$ & $0.8473$ & $1.6382$&$1.9334$\\
$300$ & $1.2291$ & $2.4377$&$1.9833$\\
$400$ & $1.6639$ & $3.2492$&$1.9528$\\
$500$ & $2.1976$ & $4.0544$&$1.8449$\\
$600$ & $2.6654$ & $5.1508$&$1.9325$\\
$700$ & $2.9765$ & $5.8691$&$1.9718$\\
$800$ & $3.3367$ & $6.5671$&$1.9681$\\
$900$ & $3.8411$ & $7.6630$&$1.9950$\\
$1000$ & $4.2435$ & $8.4006$&$1.9796$\\
\bottomrule
\end{tabular}
\end{table}

\begin{figure}
\centering
\includegraphics[width=0.8\textwidth]{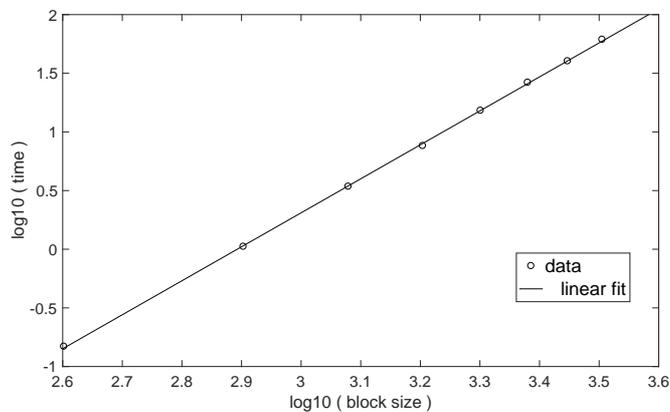}
\caption{Log-log plot of running times versus block size $m$ for Example \ref{ex_time_blocks}. The equation of the linear fit is $y=2.9x-8.4$, which confirms that the complexity of the algorithm grows as $m^3$.}\label{fig_time_blocks}
\end{figure}

\begin{example}\label{ex_time_blocks}
In this example we study the complexity of our algorithm w.r.t. block size (that is, the parameter $m$ at the end of Section \ref{sec_algo}). We choose $N=2$, $\ell=2$, $r_L=r_U=1$, with random quasiseparable generators, and focus on large values of $m$. Running times, each of them averaged over ten trials, are shown in Table \ref{table_time_blocks}. A log-log plot  is given in Figure \ref{fig_time_blocks}, together with a linear fit, which shows that the experimental growth  of these running times is consistent with theoretical complexity estimates.

\end{example}

\begin{table}
\caption{Running times in seconds for Example \ref{ex_time_blocks}.\label{table_time_blocks}}
\centering
\begin{tabular}{|c|c|}
\toprule
 block size $m$&running time (sec) \\ 
 \midrule
$400$ & $0.1481$ \\
$800$ & $1.0527$ \\
$1200$ & $3.4238$ \\
$1600$ & $7.6118$ \\
$2000$ & $15.2020$ \\
$2400$ & $26.3695$ \\
$2800$ & $40.1458$ \\
$3200$ & $61.2488$ \\
\bottomrule
\end{tabular}
\end{table}

\section{Conclusion}
In this paper we have presented an effective algorithm based on the QR decomposition for solving
a possibly large number of shifted  quasiseparable systems.  Two main motivations are the computation of a
meromorphic function of a quasiseparable matrix and the solution of linear matrix equations with quasiseparable
matrix coefficients.  The experiments performed suggest  that our algorithm has  good numerical properties. 

Future  work includes the analysis of methods based on the Mittag-Leffler's theorem for computing
meromorphic functions of quasiseparable matrices.  Approximate expansions of Mittag-Leffler type can be obtained by using  the Carath\'eodory-Fej\'er approximation theory (see \cite{garrappa2011use}). The application of these  techniques  for  computing quasiseparable matrix functions is an ongoing  research.

\bibliographystyle{plain}
\bibliography{matrixbib}

\end{document}